\documentclass[10pt,journal]{IEEEtran}
\usepackage[font=small, labelsep=space]{caption}
\usepackage{verbatim}
\usepackage{graphicx}
\usepackage{multirow}
\usepackage{amssymb}
\usepackage{amsmath}
\usepackage{amsthm}
\usepackage{amsfonts}
\usepackage{paralist}
\usepackage{stfloats}
\usepackage{url}
\usepackage{epstopdf}
\usepackage{float}
\usepackage{color}
\usepackage{cite}
\usepackage{paralist}
\usepackage[noend]{algpseudocode}
\usepackage{algorithmicx}
\usepackage{algorithm}
\usepackage{mathrsfs}
\usepackage{amsmath}
\usepackage{hyperref}
\hypersetup{hypertex=true,
            linkcolor=red,
            anchorcolor=red,
            citecolor=red}
\usepackage{subfigure}
\usepackage{amsmath}
\allowdisplaybreaks[4]
\let\itemize\compactitem
\let\enditemize\endcompactitem
\let\enumerate\compactenum
\let\endenumerate\endcompactenum

\begin{document}
\title{\LARGE{UAV-Mounted Movable Antenna: Joint Optimization of UAV Placement and Antenna Configuration}}
\author{Xiao-Wei~Tang, Yunmei Shi, Yi Huang, and Qingqing Wu
\thanks{Xiao-Wei Tang, Yunmei Shi, Yi Huang ({\{xwtang, ymshi, and huangyi718b\}@tongji.edu.cn}) are with the Department of Information and Communication Engineering, Tongji University, Shanghai, China.}
\thanks{Qingqing Wu (qingqingwu@sjtu.edu.cn) is with the Department of Electronic Engineering, Shanghai Jiao Tong University, Shanghai 200240, China.}
}
\maketitle
\begin{abstract}
Recently, movable antennas (MAs) have garnered immense attention due to their capability to favorably alter channel conditions through agile movement. In this letter, we delve into a spectrum sharing system enabled by unmanned aerial vehicle (UAV) mounted MAs, thereby introducing a new degree of freedom vertically alongside the horizontal local mobility for MAs. Our objective is to maximize the minimum beamforming gain for secondary users (SUs) while ensuring that interference to the primary users (PUs) remains below a predefined threshold, which necessitates a joint optimization involving the UAV's height, the antenna weight vector (AWV), and the antenna position vector (APV). However, the formulated optimization problem is non-convex and challenging to solve optimally. To tackle this issue, we propose an alternating optimization algorithm that optimizes the UAV's height, APV and AWV in an iterative manner, thus yielding  a near-optimal solution. Numerical results demonstrate the superiority of the proposed scheme as well as its ability to deliver full beamforming gain to SUs with reduced computational complexity.
\end{abstract}
\begin{IEEEkeywords}
MA, UAV, beamforming, spectrum sharing.
\end{IEEEkeywords}
\vspace{-3mm}
\section{Introduction}
Recently, moveable antenna (MA), also referred to as fluid antenna, has emerged as a pivotal technology endowed with the capability to be dynamically repositioned in response to changing environmental conditions or communication demands via the flexible movement of antennas \cite{Wang}. This versatility bestows upon MA several distinct advantages, including enhanced sensing accuracy \cite{Ma3}, improved system capacity \cite{Ma2}, and reduced network interference \cite{Wang2}. Consequently, the investigation into MA can thoroughly unveil the full potential of future communication systems, particularly in dynamic and unpredictable environments.

Preliminary studies have demonstrated the superiority of MAs from various perspectives. MA-enhanced multiuser communication was investigated in \cite{Zhu2}, aiming to minimize the total transmit power of users via jointly optimizing the positions of MAs, the transmit power of each user and the receive combining matrix of base station, while adhering to a minimum-achievable-rate requirement for each user. The work in \cite{Zhu3} developed a field-response model for MA-based multi-path channel by leveraging the amplitude, phase, and angle of arrival/angle of departure, based on which the achievable maximum channel gain could be largely improved compared to the conventional fixed-position antennas (FPAs) case. MA-assisted spectrum sharing was studied in \cite{Wei}, where the beamforming design and MA positions were jointly optimized to maximize the received signal power at a secondary user (SU) subject to constraints on its imposed co-channel interference power with multiple primary users (PUs). MA-based multi-beamforming \cite{Ma} focused on maximizing the minimum beamforming gain over multiple desired directions through the collaborative optimization of the antenna position vector (APV) and antenna weight vector (AWV), taking into account the restrictions on the maximum interference power over undesired directions. Despite claiming that full beamforming gains can be achieved over all desired directions, our investigation finds that there is a discrepancy between the actual and asserted beamforming gains, which may arise from neglecting the non-negativity requirement of beamforming gain when loosening the constraint associated with the APV.

The aforementioned studies mainly exploit the local movement of MAs to create favorable channel conditions. Nevertheless, by mounting the MA array onto a UAV, we introduce an additional degree of freedom beyond UAV's mobility, which allows for dynamic adjustments in the relative positions between the MAs and the users \cite{Tang}. Inspired by \cite{Ma}, we delve into multi-beamforming using a UAV-mounted MA (UMA) array to facilitate spectrum sharing services for multiple ground-based SUs. In contrast to conventional UAV-mounted base station \cite{Wu}, which prioritizes reducing {\bf{distance-dependent}} pathloss, the investigated UMA system aims to enhance {\bf{phase-sensitive}} beamforming gain by strategic UAV position adjustments. Specifically, by jointly optimizing the UAV placement and antenna configuration, we aim to maximize the minimum beamforming gain for SUs while ensuring their maximum interference to PUs. The formulated optimization problem exhibits non-convexity with respect to (w.r.t.) the UAV height, the APV, and the AWV, posing a significant challenge for solving it. To overcome this, we devise a low-complexity alternating algorithm that iteratively refines one of these variables while fixing the others. Numerical results demonstrate that the proposed algorithm enables SUs to harness the full potential of beamforming gain while effectively mitigating interference towards PUs concurrently.

\emph{Notations}: $(\cdot)$, $(\cdot)^{T}$, and $(\cdot)^{H}$ are used to denote conjugate, transpose, and conjugate transpose, respectively. The real part of vector $\boldsymbol{a}$ is denoted by $\rm{Re}\{\boldsymbol{a}\}$. $\rm{Tr}(\boldsymbol{A})$ denotes the trace of matrix $\boldsymbol{A}$. $f'(x)$ and $f''(x)$ represent the first-order and second-order derivatives of $f(x)$, respectively.

\vspace{-2mm}
\section{System Model and Problem Formulation}
\begin{figure}[htbp!]
\setlength{\abovedisplayskip}{-3mm}
\setlength{\belowcaptionskip}{-7mm}
\centering
\includegraphics[width=0.23\textwidth]{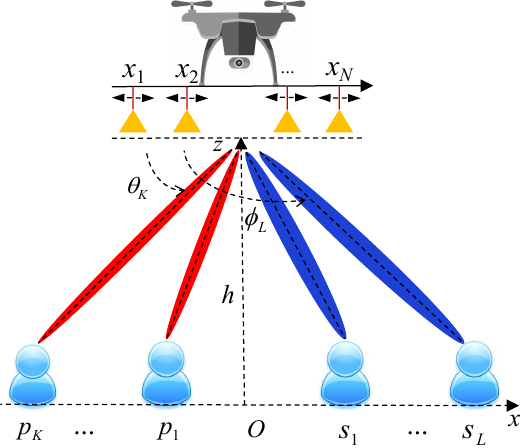}
\caption{The spectrum sharing system enabled by UAV-mounted MAs.}
\label{F1}
\end{figure}
As shown in Fig. {\ref{F1}}, we consider a spectrum sharing system enabled by UAV-mounted MAs, which consists of $K$ PUs and $L$ SUs. The locations of PUs and SUs are fixed and denoted by $\boldsymbol{p} \buildrel \Delta \over = [{p_1}, \ldots ,{p_K}]^T$ and $\boldsymbol{s} \buildrel \Delta \over = [{s_1}, \ldots, {s_L}]^T$, respectively. We assume that there are $N$ MAs, which can be flexibly deposited along the $x$-axis within a line region of length $D$. Let $\mathcal{N} \buildrel \Delta \over = \{1,\ldots, N\}$ denote the set of the MAs and ${x_n} \in [-\frac{D}{2},\frac{D}{2}]$ denote the $n$-th MA's position. Then, the APV can be denoted by ${\boldsymbol{x}} \buildrel \Delta \over = {[{x_1},{x_2}, \ldots,{x_N}]^T}$. The UAV is hovering over the origin point of $x$-axis, i.e., $x = 0$, and its height $h$ can be dynamically adjusted as needed along the $z$-axis. As a result, the steering angles over the $k$-th PU and the $l$-th SU, denoted by $\theta_k \in [0,\pi], \forall k$ and $\phi_l\in [0,\pi], \forall l$, can be respectively expressed as
\vspace{-2mm}
\begin{equation}\label{E1}
\!\!{\theta _k} = \arccos {\frac{{{-p_k}}}{\sqrt{(h^2+p_k^2)}}}, {\rm{and}}~
{\phi _l} = \arccos {\frac{{{-s_l}}}{\sqrt{(h^2+s_l^2)}}}.
\end{equation}
Hence, the steering vector can be expressed as
\begin{subequations}\label{E2}
\begin{align}
{\boldsymbol{\alpha}} ({\boldsymbol{x}},\theta_k) = {[ {{e^{j\frac{{2\pi }}{\lambda }{x_1}\cos (\theta_k)}}, \ldots ,{e^{j\frac{{2\pi }}{\lambda }{x_N}\cos (\theta_k)}}}]^T }, \forall k, \\
{\boldsymbol{\alpha}} ({\boldsymbol{x}},\phi_l ) = {[ {{e^{j\frac{{2\pi }}{\lambda }{x_1}\cos (\phi_l )}}, \ldots ,{e^{j\frac{{2\pi }}{\lambda }{x_N}\cos (\phi_l)}}} ]^T}, \forall l,
\end{align}
\end{subequations}
where $\lambda$ is the wavelength. Let ${\boldsymbol{w}}\buildrel \Delta \over = {[{w_1},{w_2}, \ldots,{w_N}]^T} \in {{\mathbb{C}}^{N}}$ denote the AWV. Thus, the beamforming gain over the steering angles $\theta_k$ and $\phi_l$ can be respectively represented as
\begin{subequations}\label{E3}
\begin{align}
G({\boldsymbol{w}},{\boldsymbol{x}},\theta_k) = {\left| {{{\boldsymbol{w}}^H}{\boldsymbol{\alpha}} ({\boldsymbol{x}},\theta_k )} \right|^2},\forall k, \\
G({\boldsymbol{w}},{\boldsymbol{x}},\phi_l) = {\left| {{{\boldsymbol{w}}^H}{\boldsymbol{\alpha}} ({\boldsymbol{x}},\phi_l )} \right|^2}, \forall l.
\end{align}
\end{subequations}

In this letter, we aim to maximize the minimum beamforming gain over SUs, denoted by $\delta$, via jointly optimizing the UAV's height $h$, the APV $\boldsymbol{x}$, and the AWV $\boldsymbol{w}$, subject to constraints on the distance between adjacent MAs, the interference towards PUs, the total available power as well as the minimum hovering height. Accordingly, the problem is formulated as
\begin{subequations}
\begin{align}
(\textbf{P1})~~&\mathop {\max }\limits_{h,{\boldsymbol{w}},{\boldsymbol{x}},\delta } \;\;\delta \notag \\
{\rm{s.t.}}~~&{x_n} - \;{x_{n - 1}} \ge {D_0},\forall n \in {\mathcal{N}}\backslash 1,\label{E4a}\\
&G({\boldsymbol{w}},{\boldsymbol{x}},{\phi _l}) \ge \delta ,\forall l \in {\mathcal{L}}, \label{E4b}\\
&G({\boldsymbol{w}},{\boldsymbol{x}},{\theta _k}) \le \eta ,\forall k \in {\mathcal{K}}, \label{E4c}\\
&\|{{\boldsymbol{w}}}\|_2  \le 1, \label{E4d}\\
&- D/2 \le {x_n} \le D/2,\forall n \in {\mathcal{N}}, \label{E4f}\\
&h \ge {H_0},\label{E4g}
\end{align}
\end{subequations}
where $D_0$ is the minimum distance for each two adjacent MAs, $\eta$ is a pre-defined interference threshold towards PUs, and $H_0$ is the minimum hovering height for the UAV. Specifically, (\ref{E4a}) ensures that there is no coupling among the MAs. (\ref{E4b}) guarantees that the beamforming gain over any SU is maintained above $\delta$. Conversely, (\ref{E4c}) imposes a constraint that the interference towards any PU must not exceed a predefined threshold, i.e., $\eta$. (\ref{E4d}) specifies that the normalized power of MAs is no larger than 1. (\ref{E4f}) ensures that the MAs should be adjusted within the confined region. (\ref{E4g}) ensures that the UAV hovers above the minimum safe height. (P1) is a formidable optimization problem, primarily due to the non-convex nature of (\ref{E4a}), (\ref{E4b}) and (\ref{E4c}) w.r.t. $h$, $\boldsymbol{w}$ or $\boldsymbol{x}$. This intricacy is further compounded by the intricate interdependence among these variables and thus significantly increases the complexity for solving (P1).

\vspace{-3mm}
\section{Proposed Algorithm}
In this section, we divide (P1) into three subproblems and  solve them iteratively in a sequential manner, where each subproblem is dedicated to optimizing either $h$, $\boldsymbol{w}$, or $\boldsymbol{x}$.

\vspace{-4mm}
\subsection{Optimization of $h$ with Given $\boldsymbol{w}$ and $\boldsymbol{x}$}
With given $\boldsymbol{w}$ and $\boldsymbol{x}$, we aim to optimize $h$ in (P1), thereby formulating the following subproblem:
\begin{align}
(\textbf{P1.2})~~&\mathop {\max }\limits_{h,\delta} \;\;\delta \notag \\
{\rm{s.t.}}~~& (\rm{\ref{E4b}}), (\rm{\ref{E4c}}), (\rm{\ref{E4g}}), \notag
\end{align}
where (\ref{E4b}) and (\ref{E4c}) are non-convex w.r.t. $h$. Hence, we relax them by adopting the successive convex approximation (SCA) technique \cite{Ma}. For ease of exposition, we denote the $n$-th element of $\boldsymbol{w}$ by ${w_n} = |{w_n}|{e^{j\angle {w_n}}}$ with amplitude $|{w_n}|$ and phase $\angle {w_n}$. Furthermore, we define ${\chi _{n,m}} \buildrel \Delta \over = \frac{{2\pi }}{\lambda }({x_n} - {x_m})$ and ${\varpi _{n,m}} \buildrel \Delta \over = \angle {w_n} - \angle {w_m}$. Thus, $G({\boldsymbol{w}},{\boldsymbol{x}},{\phi _l})$ can be further expressed as
\begin{align}\label{E5}
\!\!G({\boldsymbol{w}},{\boldsymbol{x}},{\phi _l}) = \sum\nolimits_{n = 1}^N {\sum\nolimits_{m = 1}^N {{\kappa _{n,m}}\cos ({{\hat \gamma }_{n,m,l}}(h))} } , \forall l,
\end{align}
where ${\kappa _{n,m}} \buildrel \Delta \over = |{w_n}||{w_m}|$ and ${{\hat \gamma }_{n,m,l}}(h) \buildrel \Delta \over = {\chi _{n,m}}\cos ({\phi _l}) - {\varpi _{n,m}}$. Since $G({\boldsymbol{w}},{\boldsymbol{x}},{\phi _l})$ is neither convex or concave w.r.t. $h$, we construct a surrogate function to locally approximate it based on the second-order Taylor expansion. Specifically, for a given point $\ell_0 \in \mathbb{R}$, the second-order Taylor expansion of $\cos(f(\ell))$ can be expressed as
\begin{equation}\label{E6}
\begin{array}{l}
\!\!\!\!\!\!\cos \left(f(\ell)\right) \buildrel\textstyle.\over= \cos \left(f({\ell_0})\right) - \sin (f({\ell_0})){f'}({\ell_0})(\ell - {\ell_0})\\
\!\!\!\!\!\!\!\!{-}\frac{1}{2}\left(\cos\left(f({\ell_0})\right){\left(f'({\ell_0})\right)^2} {+}\sin\left(f({\ell_0})\right)f''({\ell_0})\right){(\ell{-}{\ell_0})^2}.
\end{array}
\end{equation}
Since $(\ell{-}\ell_0)^2{\ge}0 $ and $\cos(f({\ell_0})){(f'({\ell_0}))^2}{+}\sin (f({\ell_0}))f''({\ell_0})$$\le\sqrt{(f'({\ell_0}))^4+f''({\ell_0})^2}$ according to the Cauchy-Schwartz inequality, we can construct the concave surrogate function $\hat \hbar (\ell|\ell_0)$ to approximate $\cos(f(\ell))$ as
\begin{align}\label{E7}
&\!\!\!\!\!\cos (f(\ell)) \ge \hat \hbar (\ell|{\ell_0}) \notag \\
&\!\!\!\!\!\!\buildrel \Delta \over {=}\cos (f({\ell_0})) {-} \sin (f({\ell_0})){f'}({\ell_0})(\ell {-} {\ell_0}){-} \frac{1}{2}\hat \psi ({\ell_0}){(\ell {-} {\ell_0})^2},
\end{align}
where $\hat \psi ({\ell_0}) \buildrel \Delta \over = \sqrt {{{(f'({\ell_0}))}^4}{\rm{ + }}{{(f''({\ell_0}))}^{\rm{2}}}}$. Then, for the given ${h^i}$ in the $i$-th iteration of SCA, by letting $f(\ell) \leftarrow {{\hat \gamma }_{n,m,l}}(h)$ and $f(\ell_0) \leftarrow {{\hat \gamma }_{n,m,l}}(h^i)$, we can obtain $f'(\ell_0) \leftarrow \hat \gamma _{n,m,l}'({h^i}) \buildrel \Delta \over = \frac{{{\chi _{n,m}}{s_l}{h^i}}}{{{{((h^i)^2 + s_l^2)}^{\frac{3}{2}}}}}$ and $f''(\ell_0) \leftarrow \hat \gamma _{n,m,l}''({h^i}) \buildrel \Delta \over = \frac{{{\chi _{n,m}}{s_l}(s_l^2 - 2(h^i) ^2)}}{{{{((h^i)^2 + s_l^2)}^{\frac{5}{2}}}}}$. As a result, the surrogate function that provides a global lower-bound for $G({\boldsymbol{w}},{\boldsymbol{x}},{\phi _l})$ can be constructed as
\begin{align}{\label{E8}}
G({\boldsymbol{w}},{\boldsymbol{x}},{\phi _l})\ge &\sum\nolimits_{n = 1}^N {\sum\nolimits_{m = 1}^N {{\kappa _{n,m}}\hat \hbar ({{\hat \gamma }_{n,m,l}}(h)|{{\hat \gamma }_{n,m,l}}(h^i))} } \notag \\
 \buildrel \Delta \over = &\hat a_lh^2 + {{\hat b}_l}{h} + {{\hat c}_l}, \forall l,
\end{align}
where $\hat a_l, {{\hat b}_l}, {{\hat c}_l} \in \mathbb{R}, \forall l$ are given by
\begin{align}
\!\!\!&{{\hat a}_l}=- \frac{{\rm{1}}}{{\rm{2}}}\sum\nolimits_{n = 1}^N {\sum\nolimits_{m = 1}^N {{\kappa _{n,m}}{{\hat \psi }_{n,m,l}}({h^i})} },\notag \\
\!\!\!&{{\hat b}_l}=-\sum\nolimits_{n = 1}^N\!{\sum\nolimits_{m = 1}^N {{\kappa _{n,m}}[\hat \beta_{n,m,l}(h^i)} }{-} {h^i}{{\hat \psi }_{n,m,l}}({h^i})], \notag \\
\!\!\!&{{\hat c}_l}=\sum\nolimits_{n = 1}^N {\sum\nolimits_{m = 1}^N {{\kappa _{n,m}}[\cos ({{\hat \gamma }_{n,m,l}}({h^i}))} } \notag \\
&~~~~~+ \hat \beta_{n,m,l}(h^i){h^i} - \frac{{\rm{1}}}{{\rm{2}}}{{\hat \psi }_{n,m,l}}({h^i}){({h^i})^2}], \notag
\end{align}
with ${\hat \psi }_{n,m,l}({h^i}) \buildrel \Delta \over = \sqrt {\hat \gamma _{n,m,l}'{{({h^i})}^{\rm{4}}}{\rm{ + }}\hat \gamma _{n,m,l}{''}{{({h^i})}^{\rm{2}}}}$ and $\hat \beta_{n,m,l}(h^i) \buildrel \Delta \over = \sin ({\hat \gamma }_{n,m,l}({h^i}))\hat \gamma _{n,m,l}'({h^i})$.

Additionally, since (\ref{E4c}) has a similar structure as (\ref{E4b}), we can relax it by modifying the procedure of constructing the relaxed form of (\ref{E4b}). Specifically, since $(\ell-\ell_0)^2\ge 0 $ and $\cos (f({\ell_0})){(f'({\ell_0}))^2} + \sin (f({\ell_0}))f''({\ell_0}) \ge {-}\sqrt{(f'({\ell_0}))^4{+}f''({\ell_0})^2}$ according to the Cauchy-Schwartz inequality, a global upper-bound for $G({\boldsymbol{w}},{\boldsymbol{x}},{\theta _k})$ can be approximated as
\begin{align}{\label{E9}}
G({\boldsymbol{w}},{\boldsymbol{x}},{\theta _k})\le &\sum\nolimits_{n = 1}^N {\sum\nolimits_{m = 1}^N {{\kappa _{n,m}}\bar \hbar ({{\bar \gamma }_{n,m,k}}(h)|{{\bar \gamma }_{n,m,k}}(h^i))} } \notag \\
 \buildrel \Delta \over = &\bar a_kh^2 + {{\bar b}_k}h + {{\bar c}_k}, \forall k,
\end{align}
where $\bar a_k, {{\bar b}_k}, {{\bar c}_k} \in \mathbb{R}, \forall k$ are given by
\begin{align}
\!\!\!&\bar a_k= \frac{{\rm{1}}}{{\rm{2}}}\sum\nolimits_{n = 1}^N {\sum\nolimits_{m = 1}^N {{\kappa _{n,m}}{{\bar \psi }_{n,m,k}}({h^i})} },\notag \\
\!\!\!&{{\bar b}_k} ={-}\!\!\sum\nolimits_{n = 1}^N \!{\sum\nolimits_{m = 1}^N \!{{\kappa _{n,m}}[\bar \beta_{n,m,k}(h^i) {-}{h^i}{{\bar \psi }_{n,m,k}}({h^i})]} }, \notag \\
\!\!\!&{{\bar c}_k} = \sum\nolimits_{n = 1}^N \!{\sum\nolimits_{m = 1}^N\!{{\kappa _{n,m}}[\cos ({{\bar \gamma }_{n,m,k}}({h^i})){+}\bar \beta_{n,m,k}(h^i)} } \notag \\
 &~~~~~~ + \frac{{\rm{1}}}{{\rm{2}}}{{\bar \psi }_{n,m,k}}({h^i}){({h^i})^2}], \notag
\end{align}
with ${{\bar \psi }_{n,m,k}}({h^i}){\rm{ \buildrel \Delta \over = }}\sqrt {\bar \gamma _{n,m,k}'{{({h^i})}^{\rm{4}}}{\rm{ + }}\bar \gamma_{n,m,k}''{{({h^i})}^{\rm{2}}}}$ and $\bar \beta_{n,m,k}(h^i) \buildrel \Delta \over = \sin ({\bar \gamma }_{n,m,k}({h^i}))\bar \gamma _{n,m,k}'({h^i})$.

Therefore, in the $i$-th iteration of SCA, $h$ can be optimized by solving the following optimization problem:
\begin{subequations}
\begin{align}
(\textbf{P1.2.1})~~&\mathop {\max }\limits_{h, \delta} \;\;\delta \notag \\
{\rm{s.t.}}~~&\hat a_lh^2 + {{\hat b}_l}h + {{\hat c}_l} \ge \delta ,\forall l, \label{E10a} \\
&\bar a_k h^2 + {{\bar b}_k}h + {{\bar c}_k} \le \eta ,\forall k,\label{E10b} \\
& (\rm\ref{E4g}). \notag
\end{align}
\end{subequations}
Since (\ref{E10a}) and (\ref{E10b}) are convex quadratic constraints and (\ref{E4g}) is a
linear constraint w.r.t. $h$, (P1.2.1) is a convex problem and can be efficiently solved by existing solvers, e.g., CVX.

\vspace{-2mm}
\subsection{Optimization of $\boldsymbol{w}$ with Given $h$ and $\boldsymbol{x}$}
With given $h$ and $\boldsymbol{x}$, we aim to optimize $\boldsymbol{w}$ in (P1), which leads to the following subproblem:
\begin{align}
(\textbf{P1.3})~~&\mathop {\max }\limits_{{\boldsymbol{w}},\delta} \;\;\delta \notag \\
{\rm{s.t.}}~~&(\rm{\ref{E4b}}), (\rm{\ref{E4c}}),(\rm{\ref{E4d}}), \notag
\end{align}
where (\ref{E4b}) is non-convex w.r.t. $\boldsymbol{w}$. Thus, we adopt the SCA technique to relax it. Specifically, for the given $\boldsymbol{w}^i \in \mathbb{C}^N$ in the $i$-th iteration of SCA, since $G({\boldsymbol{w}},{\boldsymbol{x}},{\phi _l})$ is convex w.r.t. $\boldsymbol{w}$, we can construct the following linear surrogate function $\bar G({\boldsymbol{w}},{\boldsymbol{x}},{\phi _l}|{{\boldsymbol{w}}^i})$ to globally approximate $ G({\boldsymbol{w}},{\boldsymbol{x}},{\phi _l})$ by applying the first-order Taylor expansion at ${{\boldsymbol{w}}^i}$:
\begin{align}\label{E11}
&G({\boldsymbol{w}},{\boldsymbol{x}},{\phi _l}) \ge \bar G({\boldsymbol{w}},{\boldsymbol{x}},{\phi _l}{\rm{|}}{{\boldsymbol{w}}^i}) \\
\buildrel \Delta \over = &2{\mathop{\rm Re}\nolimits} \{ {({{\boldsymbol{w}}^i})^H}{{\boldsymbol{\alpha }}}({\boldsymbol{x}}, {\phi _l}){{\boldsymbol{\alpha }}}{({\boldsymbol{x}},{\phi _l})^H}{\boldsymbol{w}}\}  {-} G({{\boldsymbol{w}}^i},{\boldsymbol{x}},{\phi _l}),\forall l. \notag
\end{align}
Hence, for the given ${\boldsymbol{w}}^i \in {\mathbb{C}}^N$ in the $i$-th iteration of SCA, ${\boldsymbol{w}}$ can be optimized by solving the following problem:
\begin{subequations}
\begin{align}
(\textbf{P1.3.1})~~&\mathop {\max }\limits_{{\boldsymbol{w}},\delta } \;\;\delta \notag \\
{\rm{s.t.}}~~&\bar G({\boldsymbol{w}},{\boldsymbol{x}},{\phi _l}{\rm{|}}{{\boldsymbol{w}}^i}) \ge \delta ,\forall l, \label{E12a}\\
&(\rm{\ref{E4c}}), (\rm{\ref{E4d}}), \notag
\end{align}
\end{subequations}
where (\ref{E12a}) is a linear constraint and (\ref{E4c}) and (\ref{E4d}) are convex quadratic constraints w.r.t. $\boldsymbol{w}$. Thus, (P1.3.1) is a convex problem, which can be solved via existing solvers, e.g., CVX.

\vspace{-2mm}
\subsection{Optimization of $\boldsymbol{x}$ with Given $h$ and $\boldsymbol{w}$}
With given $h$ and $\boldsymbol{w}$, we aim to optimize $\boldsymbol{x}$ in (P1), thus yielding the following subproblem:
\begin{align}
(\textbf{P1.4})~~&\mathop {\max }\limits_{{\boldsymbol{x}},\delta } \;\;\delta \notag \\
{\rm{s.t.}}~~&(\rm{\ref{E4a}}),(\rm{\ref{E4b}}),(\rm{\ref{E4c}}),(\rm{\ref{E4f}}), \notag
\end{align}
where (\ref{E4b}) and (\ref{E4c}) are non-convex constraints w.r.t. $\boldsymbol{x}$. Hence, we relax them by adopting the SCA technique. For ease of exposition, we define ${\vartheta _l} \buildrel \Delta \over = \frac{{2\pi }}{\lambda }\cos ({\phi _l})$. Therefore, $G({\boldsymbol{w}},{\boldsymbol{x}}, {\phi _l})$ can be further expressed as
\begin{align}\label{E13}
\!\!\!\!\!\!\!\!G({\boldsymbol{w}},{\boldsymbol{x}},{\phi _l}){=}\sum\nolimits_{n = 1}^N {\sum\nolimits_{m = 1}^N {{\kappa _{n,m}}\cos ({f_l}({x_n},{x_m}))} }, \forall l,
\end{align}
where ${f_l}({x_n},{x_m}) \buildrel \Delta \over = {\vartheta _l}({x_n} - {x_m}) - (\angle {w_n} - \angle {w_m})$.

Since $G({\boldsymbol{w}},{\boldsymbol{x}},{\phi _l})$ is neither convex or concave w.r.t. ${f_l}({x_n},{x_m})$, we can construct a surrogate function to locally approximate it based on the second-order Taylor expansion. Specifically, for a given $\ell_0 \in \mathbb{R}$, the second-order Taylor expansion of $\cos(\ell)$ can be expressed as
\begin{equation}\label{E14}
\!\!\!\cos(\ell)\buildrel\textstyle.\over= \cos ({\ell_0}){-} \sin ({\ell_0})(\ell {-} {\ell_0}) {-} \frac{1}{2}\cos ({\ell_0}){(\ell {-} {\ell_0})^2}.
\end{equation}
Since $\cos ({\ell_0}) \le 1$ and $(\ell-\ell_0)^2\ge 0 $, we can construct the concave surrogate function $\hat \rho (\ell|\ell_0)$ to approximate $\cos(\ell)$ as
\begin{equation}\label{E15}
\!\!\!\cos (\ell) \ge \hat \rho (\ell|\ell_0) \buildrel \Delta \over = \cos ({\ell_0}){-} \sin ({\ell_0})(\ell{-} {\ell_0}) {-} \frac{1}{2}{(\ell {-} {\ell_0})^2}.
\end{equation}

Since $G({\boldsymbol{w}},{\boldsymbol{x}},{\phi _l})$ is neither convex or concave w.r.t. $\boldsymbol{x}$, we can construct a convex surrogate function to locally approximate it based on the second-order Taylor expansion similar to Section III-A. Then, for the given ${{\boldsymbol{x}}^i} \buildrel \Delta \over = {[x_1^i,x_2^i,...,x_N^i]^T}$ in the $i$-th iteration of SCA, by letting $\ell \leftarrow {f_l}({x_n},{x_m})$ and ${\ell_0} \leftarrow {f_l}(x_n^i,x_m^i)$ in $\hat \rho (\ell|\ell_0)$ as shown in (\ref{E15}), the surrogate function that provides a global lower-bound for $G({\boldsymbol{w}},{\boldsymbol{x}},{\phi _l})$ can be constructed as
\begin{align}\label{E16}
G({\boldsymbol{w}},{\boldsymbol{x}},{\phi _l}) \ge &\sum\nolimits_{n = 1}^N {\sum\nolimits_{m = 1}^N {{\kappa _{n,m}}\hat\rho({f_l}({x_n},{x_m})|{f_l}(x_n^i,x_m^i))} } \notag \\
\buildrel \Delta \over = & \frac{1}{2}{{\boldsymbol{x}}^T}{\boldsymbol{A}_l}{\boldsymbol{x}} + {\boldsymbol{b}}_l^T{\boldsymbol{x}} + {c_l},\forall l,
\end{align}
where $\boldsymbol{A}_l \in \mathbb{R}^{N\times N}$, ${\boldsymbol{b}}_l \in \mathbb{R}^N$, and $c_l \in \mathbb{R}$ are given by
\begin{align}
&{{\boldsymbol{A}}_l} \buildrel \Delta \over =  - 2{\vartheta_l}^2(\gamma { \rm{diag}}({\boldsymbol{\bar w}}) - {\boldsymbol{\bar w}}{{{\boldsymbol{\bar w}}}^T}), \forall l, \notag \\
&{{\boldsymbol{b}}_l}[n] \buildrel \Delta \over = 2\vartheta _l^2{\sum\nolimits_{m = 1}^N {{\kappa _{n,m}}(x_n^i - x_m^i)} } \notag \\
&\;\;\;\;\;\;\; - 2{\vartheta _l} {\sum\nolimits_{m = 1}^N {{\kappa _{n,m}}\sin ({f_l}(x_n^i,x_m^i)} }, \forall l,\notag \\
&{c_l} \buildrel \Delta \over = \sum\nolimits_{n = 1}^N {\sum\nolimits_{m = 1}^N {{\kappa _{n,m}}\cos ({f_l}(x_n^i,x_m^i))} } \notag \\
&\;\;\;\;\;\; + {\vartheta _l}\sum\nolimits_{n = 1}^N {\sum\nolimits_{m = 1}^N {{\kappa _{n,m}}\sin ({f_l}(x_n^i,x_m^i))(x_n^i - x_m^i)} } \notag \\
&\;\;\;\;\;\;\;\;\; - \frac{1}{2}\vartheta _l^2\sum\nolimits_{n = 1}^N {\sum\nolimits_{m = 1}^N {{\kappa _{n,m}}{{(x_n^i - x_m^i)}^2}} }, \forall l, \notag
\end{align}
with ${\boldsymbol{\bar w}} \buildrel \Delta \over = {[|{w_1}|,|{w_2}|,...,|{w_n}|]^T}$ and $\gamma \buildrel \Delta \over = \sum\nolimits_{n = 1}^N {|{w_n}|}$. Note that ${{\boldsymbol{A}}_l}$ can be proven to be a negative semi-definite (NSD) matrix \cite{Ma}. Thus, (\ref{E4b}) is relaxed to be convex w.r.t. $\boldsymbol{x}$, that is,
\begin{equation}\label{E17}
\frac{1}{2}{{\boldsymbol{x}}^T}{\boldsymbol{A}_l}{\boldsymbol{x}} + {\boldsymbol{b}}_l^T{\boldsymbol{x}} + {c_l} \ge \delta, \forall l.
\end{equation}

On the other hand, since (\ref{E4c}) has a similar structure as (\ref{E4b}), we can relax it by modifying the procedure of constructing the relaxed convex constraint as given in (\ref{E17}). As a result, the surrogate function that provides a global upper-bound for $G({\boldsymbol{w}},{\boldsymbol{x}},{\theta _k})$ can be constructed as
\begin{align}\label{E18}
G({\boldsymbol{w}},{\boldsymbol{x}},{\theta _k}) \le & \sum\nolimits_{n = 1}^N {\sum\nolimits_{m = 1}^N {{\kappa _{n,m}}\tilde \rho({{ f}_k}({x_n},{x_m})|{{f}_k}(x_n^i,x_m^i))} } \notag \\
\buildrel \Delta \over = &\frac{1}{2}{{\boldsymbol{x}}^T}{{{\boldsymbol{\tilde A}}}_k}{\boldsymbol{x}} + \tilde {\boldsymbol{b}}_k^T{\boldsymbol{x}} + {{\tilde c}_k}, \forall k,
\end{align}
where ${{{\bf{\tilde A}}}_k} \in \mathbb{R}^{N \times N}$, ${{\tilde {\boldsymbol{b}}}_k} \in \mathbb{R}^N$, and ${\tilde c}_k \in \mathbb{R}$ are given by
\begin{align}
&{{{\boldsymbol{\tilde A}}}_k} \buildrel \Delta \over = 2\varphi _k^2(\gamma {\rm{diag}}({\boldsymbol{\bar w}}) - {\boldsymbol{\bar w}}{{{\boldsymbol{\bar w}}}^T}), \forall k,\notag \\
&{{\tilde {\boldsymbol{b}}}_k}[n] \buildrel \Delta \over =  - 2\varphi _k^2 {\sum\nolimits_{m = 1}^N {{\kappa _{n,m}}(x_n^i - x_m^i)} } \;\;\;\; \notag \\
&\;\;\;\;\;\;\;~~ - 2{\varphi _k} {\sum\nolimits_{m = 1}^N {{\kappa _{n,m}}\sin ({{f}_k}(x_n^i,x_m^i)} }, \forall k, \notag \\
&{{\tilde c}_k} \buildrel \Delta \over = \sum\nolimits_{n = 1}^N {\sum\nolimits_{m = 1}^N {{\kappa _{n,m}}\cos ({{f}_k}(x_n^i,x_m^i))} } \notag \\
&\;\;\;\;\;\; + {\varphi _k}\sum\nolimits_{n = 1}^N {\sum\nolimits_{m = 1}^N {{\kappa _{n,m}}\sin ({{f}_k}(x_n^i,x_m^i))(x_n^i - x_m^i)} } \notag \\
&\;\;\;\;\;\;\; + \frac{1}{2}\varphi _k^2\sum\nolimits_{n = 1}^N {\sum\nolimits_{m = 1}^N {{\kappa _{n,m}}{{(x_n^i - x_m^i)}^2}} }, \forall k. \notag
\end{align}
Note that ${{{\boldsymbol{\tilde A}}}_k}$ can be rigorously proven to be a positive semi-definite matrix \cite{Ma}. Thus, (\ref{E4c}) can be relaxed as a convex constraint w.r.t. $\boldsymbol{x}$:
\begin{equation}\label{E19}
\frac{1}{2}{{\boldsymbol{x}}^T}{{\tilde {\boldsymbol{A}}}_k}{\boldsymbol{x}} + {\tilde {\boldsymbol{b}}}_k^T{\boldsymbol{x}} + {{\tilde c}_k} \le \eta, \forall k.
\end{equation}

Notice that $G({\boldsymbol{w}},{\boldsymbol{x}},{\theta _k})$ must be larger than 0. However, the relaxed form of $G({\boldsymbol{w}},{\boldsymbol{x}},{\theta _k})$ as given in (\ref{E19}) can not guarantee this requirement. Thus, the following constraint should be satisfied:
\begin{equation}\label{E20}
\frac{1}{2}{{\boldsymbol{x}}^T}{{\dot {\boldsymbol{A}}}_k}{\boldsymbol{x}} + {\dot {\boldsymbol{b}}}_k^T{\boldsymbol{x}} + {{\dot c}_k} \ge 0, \forall k.
\end{equation}
where ${{{{\dot {\boldsymbol{A}}}}}_k} \in \mathbb{R}^{N \times N}$, ${{\dot {\boldsymbol{b}}}_k} \in \mathbb{R}^N$, and ${\dot c}_k \in \mathbb{R}$ are given by
\begin{align}
&{{\dot{\boldsymbol{A}}}_k} \buildrel \Delta \over = -2\varphi _k^2(\gamma {\rm{diag}}({\boldsymbol{\bar w}}) - {\boldsymbol{\bar w}}{{{\boldsymbol{\bar w}}}^T}), \forall k,\notag \\
&{{\dot {\boldsymbol{b}}}_k}[n] \buildrel \Delta \over = 2\varphi _k^2 {\sum\nolimits_{m = 1}^N {{\kappa _{n,m}}(x_n^i - x_m^i)} } \;\;\;\;\notag \\
&\;\;\;\;\;\;\; - 2{\varphi _k} {\sum\nolimits_{m = 1}^N {{\kappa _{n,m}}\sin ({{f}_k}(x_n^i,x_m^i)} }, \forall k, \notag \\
&{{\dot c}_k} \buildrel \Delta \over = \sum\nolimits_{n = 1}^N {\sum\nolimits_{m = 1}^N {{\kappa _{n,m}}\cos ({{ f}_k}(x_n^i,x_m^i))} } \notag \\
&\;\;\;\;\;\; + {\varphi _k}\sum\nolimits_{n = 1}^N {\sum\nolimits_{m = 1}^N {{\kappa _{n,m}}\sin ({{f}_k}(x_n^i,x_m^i))(x_n^i - x_m^i)} } \notag \\
&\;\;\;\;\;\;\;\;\; - \frac{1}{2}\varphi _k^2\sum\nolimits_{n = 1}^N {\sum\nolimits_{m = 1}^N {{\kappa _{n,m}}{{(x_n^i - x_m^i)}^2}} }, \forall k. \notag
\end{align}
Note that ${{\dot {\boldsymbol{A}}}_k}$ can also be proven to be a NSD similar to ${{ {\boldsymbol{A}}}_l}$. Therefore, in the $i$-th iteration of SCA, $\boldsymbol{x}$ can be optimized by solving the following optimization problem:
\begin{align}
(\textbf{P1.4.1})~~&\mathop {\max }\limits_{{\boldsymbol{x}},\delta} \;\;\delta \notag \\
{\rm{s.t.}}~~&(\rm{\ref{E4a}}),(\rm{\ref{E4f}}), (\rm{\ref{E17}}), (\rm{\ref{E19}}),(\rm{\ref{E20}}).\notag
\end{align}
Since (\ref{E4a}) and (\ref{E4f}) are linear constraints and (\ref{E17}), (\ref{E19}) and (\ref{E20}) are convex quadratic constraints w.r.t. $\boldsymbol{x}$, (P1.4.1) is a convex problem and can be efficiently solved by CVX.

\vspace{-2mm}
\subsection{Overall Algorithm and Complexity Analysis}
The overall algorithm for solving (P1) is summarized in Algorithm 1. Let $I_h$, $I_{\boldsymbol{w}}$ and $I_{\boldsymbol{x}}$ denote the number of iterations for solving (P1.2.1), (P1.3.1), and (P1.4.1), respectively. In each iteration, $h$, $\boldsymbol{w}$ and $\boldsymbol{x}$ are alternatively optimized using the interior-point method, and thus their individual complexity can be represented as $O(I_h(L+K+1)^3\ln(1/\varsigma))$, $O(I_{\boldsymbol{w}}(L+K+1)N^{3.5}\ln(1/\varsigma))$ and $O(I_{\boldsymbol{x}}(2N+2K+L-1)N^{3.5}\ln(1/\varsigma))$, respectively, with $\varsigma$ being the pre-specified precision. Hence, the total computational complexity is $O(I(I_h(L+K+1)^3+I_{\boldsymbol{w}}(L+K+1)N^{3.5}+I_{\boldsymbol{x}}(2N+2K+L-1)N^{3.5})\rm{log}(1/\varsigma))$ with $I$ denoting the number of iterations for iteratively solving (P1.2.1), (P1.3.1), and (P1.4.1).
\vspace{-2mm}
\begin{algorithm}[htb]
\caption{The alternating optimization for UMA.}
\label{alg:Framework}
\begin{algorithmic}[1]
\State Input: $N$, $L$, $K$, $\eta$, $H_0$, $\boldsymbol{p}$, $\boldsymbol{s}$,$D$,$D_0$, $\varsigma$.
\State Set iteration index $i = 0$ and initialize $h^0$, {{$\boldsymbol{w}^0$}}, and {{$\boldsymbol{x}^0$}}.
\State {\bf{repeat}}
\State Solve (P1.2.1) and denote the optimal solution to the UAV's height as ${{{h}}^{i+1}}$.
\State Solve (P1.3.1) and denote the optimal solution to the AWV as ${{{\boldsymbol{w}}}^{i+1}}$.
\State Solve (P1.4.1) and denote the optimal solution to the APV as ${{{\boldsymbol{x}}}^{i+1}}$.
\State Update $i = i+1$.
\State {\bf{until}} the computed objective value of (P1) converges within a pre-specified precision $\varsigma>0$.
\State Output: $\delta$, $h$, $\boldsymbol{w}$, $\boldsymbol{x}$.
\end{algorithmic}
\end{algorithm}

\vspace{-4mm}
\section{Numerical Results}
In the simulation, unless otherwise specified, we set $N=8$, $K=2$, $L=2$, $\eta = 0.1$, $H_0 = 10$m, $\lambda = 0.1$m, $D_0 = \frac{\lambda}{2}$, $D = 8D_0$, and $\varsigma = 10^{-3}$. The locations of SUs and PUs are set to $\boldsymbol{s} = [-11.91,5.77]^T$m and $\boldsymbol{p} = [-56.71,17.32]^T$m. $\boldsymbol{w}^0$ and $\boldsymbol{x}^0$ can be initialized by referring to \cite{Ma}, while $h^0$ can be initialized via $h^0 = \frac{\sum\nolimits_{n = 1}^N\sum\nolimits_{m = 1}^N\sum\nolimits_{l = 1}^L h_{n,m,l}^0}{LN^2}$. Specifically, $h_{n,m,l}^0$ can be obtained via solving the following equality:
\begin{figure}[htbp!]
\setlength{\abovedisplayskip}{-6mm}
\setlength{\belowcaptionskip}{-7mm}
\centering
\includegraphics[width=0.27\textwidth]{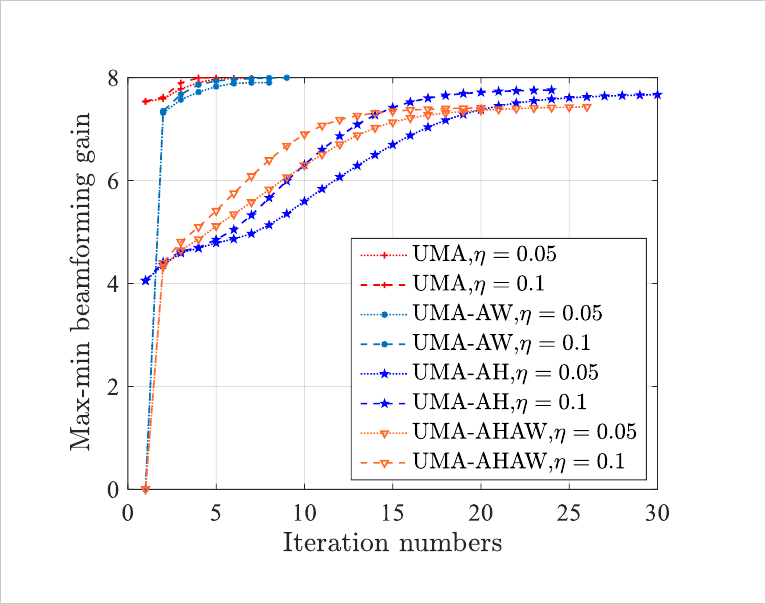}
\caption{Max-min beamforming gain versus iteration numbers.}
\label{F2}
\end{figure}
\begin{equation}
\tan (\frac{{ - {\chi _{n,m}}{s_l}}}{{\sqrt {z_{n,m,l}} }} - {\varpi _{n,m}}) = \frac{{(s_l^2 - 2{{({h_{n,m,l}^0})}^2}){{(z_{n,m,l})}^{\frac{1}{2}}}}}{{{\chi _{n,m}}{s_l}{{({h_{n,m,l}^0})}^2}}},
\end{equation}
where $z_{n,m,l} \buildrel \Delta \over = h_{n,m,l}^0+s_l^2$. In addition, we compare the proposed scheme with 3 benchmarks named UMA-AH, FPA, MA, whose details are given as: 1) \textbf{UMA-AH}: The initial UAV height of the proposed algorithm is set arbitrarily, i.e., $h^0 = H_0$; 2) \textbf{FPA}: The height of the UAV and the positions of the MAs are fixed; 3) \textbf{MA}\cite{Ma}: The height of the UAV is fixed while the AWV and APV are alternatively optimized.

In Fig. \ref{F2}, we show the max-min beamforming gain versus iteration numbers for the proposed scheme and its variants. Here, `AH/AW' is abbreviated for {\underline{a}}rbitrary {\underline{h}}eight/A{\underline{W}}V. Notably, the proposed UMA scheme initiates from an exceedingly close-to-optimal starting point and exhibits swift convergence towards the full beamforming gain across different $\eta$ values. For example, the proposed UMA scheme only takes about $\frac{1}{5}$ iteration numbers of the UMA-AWAH scheme under $\eta = 0.1$. The UMA-AW approach experiences a marginal decrement in max-min beamforming gain under $\eta = 0.05$ and a slightly prolonged convergence period, which, however, still surpasses the UMA-AH and UMA-AHAW methods in terms of convergence speed and overall system performance. Consequently, we deduce from Fig. \ref{F2} that the proposed height initialization technique (e.g., UMA and UMA-AW) contributes significantly to enhancing convergence speed and elevating the max-min beamforming gain.

\begin{figure}[htbp!]
\setlength{\abovedisplayskip}{-2mm}
\setlength{\belowcaptionskip}{-3mm}
\centering
\includegraphics[width=0.26\textwidth]{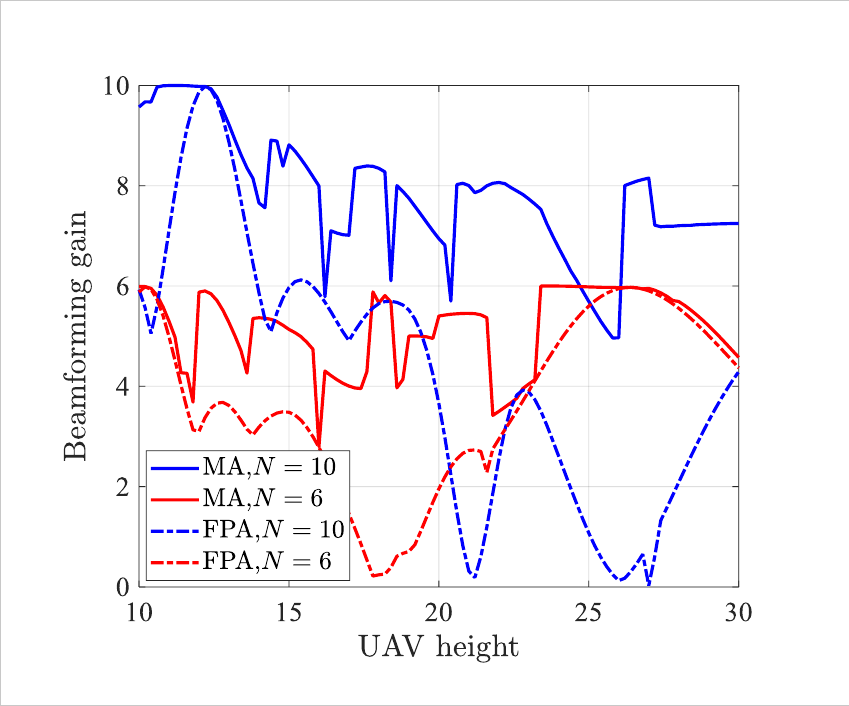}
\caption{Beamforming gain versus UAV height.}
\label{F3}
\end{figure}

In Fig. \ref{F3}, we demonstrate the beamforming gain for the 1-st SU versus UAV height with $N=6$ and $N=10$ for the FPA and MA benchmarks. The beamforming gain fluctuates sharply with the UAV height in the FPA scheme, which, however, is much more gentle in the MA's case. This is because the varying UAV height changes the relative positions between the UAV and SU, thus affecting the steering vector as defined in (\ref{E2}). It is also observed from Fig. \ref{F3} that the MA scheme can achieve the full beamforming gain for the 1-st SU within a specific height range, e.g., $h \in [10.5, 12.5]$ for $N = 10$, while the FPA scheme can only achieve the beamforming gain at a certain point, e.g., $h = 12.5$ for $N = 10$. This indicates that by adjusting the UAV height based on the MA scheme, it is highly likely to identify an optimal UAV position, thereby enabling all SUs to achieve the full beamforming gain.

\begin{figure}[htbp!]
\setlength{\abovedisplayskip}{10pt}
\setlength{\belowcaptionskip}{-3mm}
\centering
\includegraphics[width=0.27\textwidth]{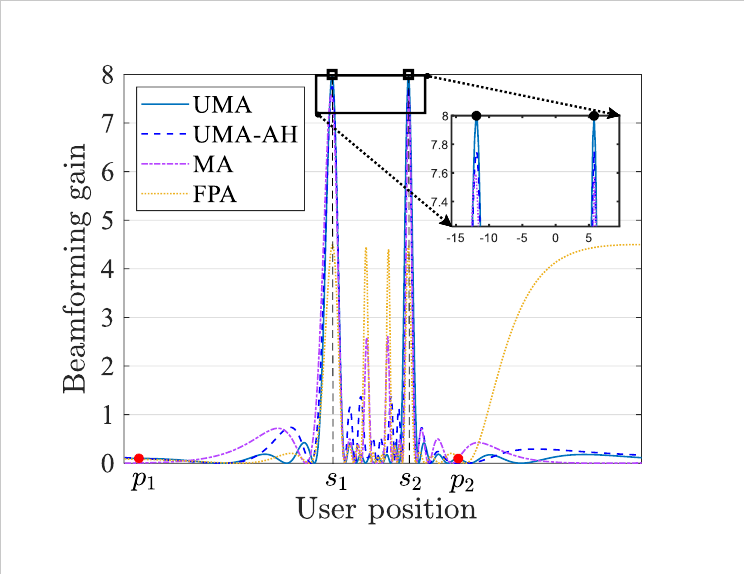}
\caption{Comparison of beam patterns with different benchmarks.}
\label{F4}
\end{figure}

Fig. \ref{F4} presents the comparison of beam patterns with different benchmarks. We can see from Fig. \ref{F4} that the beamforming gain for the two PUs can be well restrained under the pre-determined threshold for the four considered schemes. Moreover, we can see that the proposed scheme achieves the full beamforming gain for the two SUs (i.e., $G({\boldsymbol{w}},{\boldsymbol{x}},{\phi _l})=8$, $l \in \{1,2\}$) while the UMA-AH, MA, FPA counterparts achieve a beamforming gain of 7.75, 7.57, and 4.47, respectively.
This is because the proposed UMA scheme can flexibly adjust the steering vector by exploiting the new degree of freedom provided by the UAV height adjustment.

\section{Conclusions}
In this letter, we investigated a UMA system to enhance the achievable beamforming gain for SUs by exploiting the UAV mobility and local MA movement. A low-complexity alternating optimization algorithm was devised to obtain a near-optimal solution to the formulated non-convex optimization problem.  Numerical results demonstrated that the proposed UMA scheme outperformed its UMA-AH, MA and FPA counterparts, which could achieve the full beamforming gain for all SUs while mitigating the interference towards PUs simultaneously with reduced computational complexity thanks to the proposed UAV height initialization technique.

\vfill

\bibliographystyle{IEEEtran}

\vspace{-2mm}
\begin{thebibliography}{99}
\bibitem{Wang}
C. Wang et al., ``AI-empowered fluid antenna systems: Opportunities, challenges, and future directions,'' \emph{IEEE Wireless Commun.}, early access, pp. 1-8, Jul. 2024.
\bibitem{Ma3}
W. Ma, L. Zhu, and R. Zhang, ``Movable antenna enhanced wireless sensing via antenna position optimization,'' 2024, arXiv:2405.01215.
\bibitem{Ma2}
W. Ma, L. Zhu, and R. Zhang, ``MIMO capacity characterization for movable antenna systems,'' \emph{IEEE Trans. Wireless Commun.}, vol. 23, no. 4, pp. 3392-3407, Apr. 2024.
\bibitem{Wang2}
H. Wang et al., ``Movable antenna enabled interference network: Joint antenna position and beamforming design,'' \emph{IEEE Wireless Comm. Lett.}, early access, pp. 1-5, Jul. 2024.
\bibitem{Zhu2}
L. Zhu et al., ``Movable-antenna enhanced multiuser communication via antenna position optimization,'' \emph{IEEE Trans. Wireless Commun.}, vol. 23, no. 7, pp. 7214-7229, Jul. 2024.
\bibitem{Zhu3}
L. Zhu, W. Ma, and R. Zhang, ``Modeling and performance analysis for movable antenna enabled wireless communications,'' \emph{IEEE Trans. Wireless Commun.}, vol. 23, no. 6, pp. 6234-6250, Jun. 2024.
\bibitem{Wei}
X. Wei et al., ``Joint beamforming and antenna position optimization for movable antenna-assisted spectrum sharing,'' \emph{IEEE Wireless Commun. Lett.}, early access, pp. 1-5, Jul. 2024.
\bibitem{Ma}
W. Ma et al., ``Multi-beam forming with movable-antenna array,'' \emph{IEEE Commun. Lett.}, vol. 28, no. 3, pp. 697-701, Mar. 2024.
\bibitem{Tang}
X. -W. Tang et al., ``3D trajectory planning for real-time image acquisition in UAV-assisted VR,'' \emph{IEEE Trans. Wireless Commun.}, vol. 23, no. 1, pp. 16-30, Jan. 2024.
\bibitem{Wu}
Q. Wu et al., ``A comprehensive overview on 5G-and-beyond networks with UAVs: From communications to sensing and intelligence,'' \emph{IEEE J. Sel. Areas Commun.}, vol. 39, no. 10, pp. 2912-2945, Oct. 2021.
\end{thebibliography}
\end{document}